\documentclass{amsart}

\usepackage{amssymb}
\usepackage[all]{xy}
\usepackage{hyperref}

\theoremstyle{definition}

\newtheorem*{ack}{Acknowledgments}      

\newtheorem{defn-thm}[thm]{Definition--Theorem}  
\newtheorem{defn-lem}[thm]{Definition--Lemma}  

\theoremstyle{remark}

\newcommand{\z}[0]{{\mathbb Z}}

\def\loccoh#1.#2.#3.#4.{H^{#1}_{#2}(#3,#4)}

\DeclareMathAlphabet{\mathchanc}{OT1}{pzc}%
                                {m}{it}

\usepackage[all]{xy}\xyoption{dvips}

\begin{document}
\bibliographystyle{amsalpha}


\title{Arthur Byron Coble, 1878--1966}
        \author{J\'anos Koll\'ar}

        \begin{abstract}
A short essay on the life and mathematical heritage of Coble. 
          A substantially edited version will be part of the series of 
biographical memoirs of past members 
          of  the National Academy of Sciences.
        \end{abstract}

\maketitle

 Coble was born  November 3, 1878 in Dauphin County, Pennsylvania, near Harrisburg. He graduated from Pennsylvania College (now  Gettysburg College) in 1897.
 After a year of public school teaching, he studied at the Johns Hopkins University (1898--1902), earning a Ph.D.\ with  Frank Morley. The title of his dissertation was
 {\it The quartic curve as related to conics.} He taught for one year at the University of Missouri, then returned to the Johns Hopkins University as a research assistant at the Carnegie Institute, where  he was later promoted to research  associate and   to associate professor. In 1904 he visited Greifswald  and Bonn Universities in Germany  with the support of the Carnegie Institute.

 In 1918 Coble accepted a professorship at the University of Illinois
  at Urbana-Champaign and stayed there save for visiting positions at the University of Chicago (1919) and the Johns Hopkins University (1927--28). He was
  head of the department  from 1933  until his retirement at 1947.
 He moved back to 
 Dauphin County, Pennsylvania and passed away in  Harrisburg on December 8, 1966.

 Coble was very active in the  American Mathematical Society,
 served on the governing council (1911-14), vice president (1917--20),
 chair of the Chicago section (1922) and  as president (1933-34).
 He was editor of the AMS Transactions (1920--25), Proceedings of the AMS (1933-34)  and Duke Math.\ Journal (1936--38).
 Several times he served on the National Research Council   and on investigating committees of the American Association of University Professors.
 
 He was elected to the National Academy of Sciences in 1924, delivered the AMS Colloquium lectures  in 1928, and received an honorary degree
 from Gettysburg College in 1932.

Coble had  27 doctoral students, among them  7 were women, starting with  Bessie Miller (the Johns Hopkins University, 1914) and ending with 
  Janie Lapsley Bell (University of Illinois, 1943). He was among the leading advisors for women doctorates in mathematics before 1940. 
   Coble's students did not seem to have continued his work in algebraic geometry.   About half of the theses of his  over 200 descendants   are in
applied mathematics and automata  theory, the other half in mathematics education.
\medskip

{\bf Mathematical works.} 
MathSciNet lists 24 publications, Archibald's volume on the history of the AMS (1938) gives references to  43,  while the Memorial resolution of the University of Illinois Senate (1968) mentions over 60.


Coble's most important contribution is the book
{\it Algebraic geometry and theta functions}  (American Mathematical Society Colloquium Publications, Vol.\ 10, AMS, Providence, R.I., 1929).
The detailed review by  Zariski notes that
``thanks to its rich geometric content and originality of treatment, Coble's
work is a really important contribution to the theory and application of the
$\theta$-functions'' and that it is for the ``competent readers.''
({\it Algebraic Geometry and
Theta Functions,} Bull.\ Amer.\ Math.\ Soc.\ 36 (1930), no. 7,
452--454).

A more precise title would have been `Cremona transformations and theta functions,'  two major topics of  investigation at that time.
Besides giving a broad overview of the long history of these areas, the book also contains a good
description  of Coble's research, with many new results added.
The  book gives an almost complete picture of his research in this area.

The Cremona group in dimension $n$, denoted by $\operatorname{Cr}_n$, is
the group of birational automorphisms of projective $n$-space ${\mathbb P}^n$.
Since its formal definition by Cremona in 1863, it has been an object of intense study; see the books by L.~Godeaux  ({\it Les transformations birationelles du plan,} M\'em.\ des sci.\ math. vol.\ 22, Gauthier-Villars et Cie, 1927), and of 
H.P.~Hudson ({\it Cremona transformations in plane and space,} Cambridge University Press, 1927), and the numerous references in them for contemporary accounts.

A major discovery of Coble relates the plane Cremona group 
$\operatorname{Cr}_2$ to the  Coxeter-Weyl groups
of the sequence of lattices
$
A_1+A_2, A_4, D_5, E_6, E_7, E_8, E_9, E_{10}, \dots 
$
(Kodaira uses  $\tilde E_8$ for $E_9$, and, especially  for $m\geq 9$, these are also frequently denoted by $T_{3,2,m-3}$.)

Using modern terminology, let  $S_m$ be a surface obtained by blowing up $m$ points in ${\mathbb P}^2$ (in general position) and let $E_m\subset H_2(S_m, \z)$ be the orthogonal complement of the first Chern class. Poincar\'e duality gives a quadratic form, which is indefinite for $m\geq 10$. These are exactly the lattices mentioned above.

Let $W(E_m)$ denote  the
Coxeter-Weyl group of $E_m$.
Coble observes that there is a  natural (birational) action of 
$W(E_m)$ on
 the configuration space of $m$ points in ${\mathbb P}^2$.
For $m\leq 8$ these give rise to the
construction of the moduli space of Del~Pezzo surfaces of degree
$9-m$. The case $m=6$ corresponds to cubic surfaces; these have been much studied in the 19th century. Coble also gives detailed information about the $m=7$ and $m=8$ cases, the latter being especially complicated. 

For $m=10$  Coble shows that the set of  10 singular points of a rational sextic curve (sometimes called a Coble curve), is invariant under
the $W(E_{10})$-action. Equivalently, if we blow up the 10 singular points, then the automorphism group of the resulting surface is   $W(E_{10})$. A  theorem of 
S.\ Cantat and I.\ Dolgachev  ({\it Rational surfaces with a large group of automorphisms,} J. Amer. Math. Soc.
25 (2012), No. 3, 863-905), shows that these---now called Coble surfaces---are the only surfaces with this property in characteristic 0. (There is one more series in positive characteristic, found by  B.~Harbourne (1985).)

Coble also studied the higher dimensional cases, and further generalizations are in 
S.~Mukai
({\it Geometric realization of T-shaped root systems and counterexamples to Hilbert's fourteenth problem.}
Enc.\ Math.\ Sci. 132,
Springer, Berlin, 2004).


Coble surfaces also play an important role in the theory of Enriques surfaces. 
Chapter 9 of the monograph
{\it Enriques Surfaces II\footnote{\url{https://dept.math.lsa.umich.edu/~idolga/EnriquesTwo.pdf}}} by 
I.\ Dolgachev and S.\ Kond\=o  is devoted to their detailed study.

Coble also investigates a remarkable  partial inverse of the $W(E_m)$-action.
Roughly speaking, the   sub-groupoid of  $\operatorname{Cr}_2$ consisting of
Cremona transformations with $\leq m$ base points maps to  $W(E_m)$.
(As Coble notes, there is an issue here with ordered versus un-ordered sets of base points.) In some vague sense this suggests that  $\operatorname{Cr}_2$ might act on some hyperbolic space.

It turns out that $\operatorname{Cr}_2$ cannot act on a finite dimensional space, but if we take the limit, blowing up all the points, then we get a representation of  $\operatorname{Cr}_2$ on an infinite dimensional hyperbolic space.  S.\ Cantat used this to obtain a series of theorems on $\operatorname{Cr}_2$ ({\it Sur les groupes de transformations birationnelles des surfaces,} Annals of Math., Vol. 174, 2012, pp. 299–334). It is unlikely that Coble  foresaw this development,  but his work was the first to establish some connection between
$\operatorname{Cr}_2$ and hyperbolic groups.

From Coble's point of view, a key earlier result for theta functions was a theorem of Thomae (1870) 
that, in modern language, described the coordinate ring of the geometric invariant theory quotient of the configuration space of $n$ points in ${\mathbb P}^1$
using theta functions.
Another result---going back at least to Pl\"ucker (1839)---concerns the bitangents of  plane quartic curves (a topic of surprisingly many papers). Bitangents are linear objects, but they are also 
the odd theta characteristics of genus 3 curves, and   give the
generators of the Kleiman-Mori cone of degree~2 Del~Pezzo surfaces.
So a property in classical linear geometry turns out to be equivalent to
abstract properties of curves and surfaces.
Coble's AMS Colloquium lecture was devoted to the genus 4 case, where the
odd theta characteristics correspond to tritangent planes of the degree 6 canonical image of the curve in ${\mathbb P}^3$.
Much of Coble's book can be  viewed as a search for  generalizations of these examples.

Coble also worked out two subtle 
 geometric invariant theory quotients of  configuration spaces (before this notion  existed).
For 6 ordered points in   ${\mathbb P}^1$ he gets the Segre cubic, and
for 6 ordered points in   ${\mathbb P}^2$ a
 double cover of ${\mathbb P}^4$
({\it Point sets and allied Cremona groups,} Trans.\ Amer.\ Math.\ Soc.  16 (1915), 155–198).
The latter is usually called the  Coble variety.

Two remarkable types of hypersurfaces are also named after Coble, who discovered them.
One is a quartic in ${\mathbb P}^7$ whose singular locus is a Kummer 3-fold of a genus 3 curve, and the other  is a cubic in ${\mathbb P}^8$ whose singular locus is a Jacobian of a genus 2 curve.
A modern treatment is given by 
A.~Beauville ({\it The Coble hypersurfaces,} Comptes Rendus Math. 337  (2003): 189–194).

Dolgachev points out that,  for a historian of mathematics, Coble's book is frustratingly haphazard with citations. The preface of the book clearly states which sections contain new material, but  facts frequently appear in the other sections without citations and without clear indications of how the proofs should go.
The basics of  the birational  representations of the Weyl groups
were  known to S.~Kantor, and the  theory of `association' (now called
Gale transform) can be traced back to G.~Castelnuovo.
It is also not clear where the results on a
special  quartic threefold---studied by  G.~Castelnuovo (1891) and H.W.~Richmond (1902), but now named after J.-i.~Igusa (1962)---come from\footnote{As the current names of these objects indicate, Coble was not the last one to rediscover old theorems and examples.}.
\medskip

 {\bf  Coble's influence.} His book appeared at a transitional period of algebraic geometry.
 Much of the earlier work was devoted to the geometry of subvarieties of projective spaces, with emphasis on especially nice and interesting examples in low dimensions.
 Around this time van~der~Waerden,  Weil and Zariski were reworking algebraic geometry to focus on general theorems, and free it from the constraints of  projective space.

 Coble  also  moved beyond the linearity imposed by projective geometry, by looking at properties invariant under the birational automorphism group
 $\operatorname{Cr}_n$, instead of  the traditional linear automorphism group
 $\operatorname{PGL}_{n+1}$.
In hindsight, although Coble made a significant step in the right direction, it was not decisive enough.

This explains why Coble's book
was not as influential as some other  volumes in that series around the same time (for example by Lefschetz,   Morse or Stone).
 While the books by  Yu.\ I.\  Manin ({\it Cubic forms: algebra, geometry, arithmetic,} Nauka, Moscow, 1972)
and  D.\ Mumford
({\it Tata lectures on theta I--III.}  Progress in Math. vols.\ 28, 43, 97. Birkh\"auser, Boston, MA, 1983--1991)
both cite Coble in passing,  neither does justice to the  wealth of material in his book.

On the other hand,  Coble's  book  was reprinted in 1962 and again in 1980, and on MathSciNet it still gets an average of 4 citations each year.
The borrowing slips in the Princeton University library show that   copies  were in constant use.

As algebraic geometry returned to concrete questions, much of Coble's work was revisited.  His main ideas were developed with modern terminology and rigor by 
I.~Dolgachev,  whose book  ({\it Points Sets in Projective Spaces and Theta Functions,} Ast\'erisque, vol. 165, 1988, with D.~Ortland)
makes Coble's book  accessible to contemporary readers.

We can now view Coble as an important figure of american mathematics, who created a lasting legacy in algebraic gometry.


\medskip
 {\bf Sources.} The  information on Coble's life mostly comes from
 R.\ C.\ Archibald ({\it A semicentennial history of the American Mathematical Society 1888-1938,} New York, 1938, pp.233-236), and the memorial resolution of the
 University of Illinois Senate, which is reproduced by 
 A.\ Mattuck ({\it Arthur Byron Coble,}  Bull.\ Amer.\ Math.\ Soc.\ 76 (1970), 693-699). The latter also gives an appreciation of Coble's mathematical works,  
that
``serve to our present-day algebraic geometers, dwelling as they do in their Arcadias of abstraction, as a reminder of what awaits those who dare to ask specific questions about particular varieties.'' 
  Additional details can be found in 
  K.\ H.\ Parshall ({\it Arthur Byron Coble}, American National Biography 5 (Oxford, 1999), 113-114), and the connections with  Gettysburg college are treated by D.\ B.\ Glass ({\it Coble and Eisenhart: two Gettysburgians who led mathematics.} Notices Amer.\ Math.\ Soc.\ 60 (2013), no.\ 5, 558–566).

  For Coble's students, see the Mathematics Genealogy
  Project\footnote{\url{https://www.genealogy.math.ndsu.nodak.edu}}.
 More details on the women students are in the books  by J.\ Green and  J.\ LaDuke
({\it Pioneering women in American mathematics: the pre-1940 PhD's,}  
American Mathematical Society, Providence, R.I.;  London Mathematical Society, London,  2009) and
M.\ A.\ M.\ Murray ({\it Women Becoming Mathematicians: Creating a Professional Identity in Post World War II America,} MIT Press, 2000).

See also the Wikipedia
entry\footnote{\url{https://en.wikipedia.org/wiki/Arthur_Byron_Coble}} and
the History of Mathematics
Archive\footnote{\url{https://mathshistory.st-andrews.ac.uk/Biographies/Coble}}.

\begin{ack}  I thank  Robert C.\ Gunning and Margaret A.M.\ Murray for helpful comments and references.
The main expert on Coble's mathematical work is  Igor Dolgachev; his guidance and comments were invaluable in preparing this short summary.
Partial  financial support    was provided  by  the NSF under grant number
DMS-1901855.

\hfill Princeton University,  \email{kollar@math.princeton.edu}
 \end{ack}


\end{document}